\documentclass[12pt]{article}
\usepackage{amsmath,amsthm,amsfonts,amscd,amssymb} 
\usepackage[T1]{fontenc}
\usepackage[utf8]{inputenc}
\usepackage{authblk}
\usepackage{verbatim}   
\usepackage{eucal} 
\usepackage{longtable} 
\usepackage{makeidx}
\usepackage{url}

\newcommand{\field}[1]{\mathbb{#1}}
\newcommand{\califas}[1]{\mathcal{#1}}
\newcommand{\Z}{\field{Z}}
\newcommand{\Q}{\field{Q}}

\newcommand{\C}{\field{C}}

\newcommand{\F}{\field{F}}

\newcommand{\Ec}{\califas{E}}

\newcommand{\Ct}{\tilde{C}}

\newfont{\cyrr}{wncyr10}
\def\Sh{\mbox{\cyrr Sh}}

\providecommand{\keywords}[1]{\small \textbf{Keywords:} #1}
\providecommand{\subjclass}[1]{\indent {\small {\bf 2010 Mathematics Subject Classification:} #1}}
\providecommand{\email}[1]{e-mail: \texttt{#1}}
\newcommand{\thm}{\begin{theorem}}
\newcommand{\thml}[1]{\begin{theorem}\label{#1}}
\newcommand{\mht}{\end{theorem}}
\newcommand{\cnj}{\begin{conjecture}}
\newcommand{\cnjl}[1]{\begin{conjecture}\label{#1}}
\newcommand{\jnc}{\end{conjecture}}
\newcommand{\dfn}{\begin{definition}}
\newcommand{\dfnl}[1]{\begin{definition}\label{#1}}
\newcommand{\nfd}{\end{definition}}
\newcommand{\pro}{\begin{proposition}}
\newcommand{\prol}[1]{\begin{proposition}\label{#1}}
\newcommand{\orp}{\end{proposition}}
\newcommand{\crl}{\begin{corollary}}
\newcommand{\crll}[1]{\begin{corollary}\label{#1}}
\newcommand{\lrc}{\end{corollary}}
\newcommand{\lmm}{\begin{lemma}}
\newcommand{\lmml}[1]{\begin{lemma}\label{#1}}
\newcommand{\mml}{\end{lemma}}
\newcommand{\prf}{\begin{proof}}
\newcommand{\prfl}[1]{\begin{proof}\label{#1}}
\newcommand{\frp}{\end{proof}}
\newcommand{\axi}{\begin{axiom}}
\newcommand{\axil}[1]{\begin{axiom}\label{#1}}
\newcommand{\ixa}{\end{axiom}}
\newcommand{\rmk}{\begin{remark}}
\newcommand{\rmkl}[1]{\begin{remark}\label{#1}}
\newcommand{\kmr}{\end{remark}}
\newcommand{\exa}{\begin{example}}
\newcommand{\exal}[1]{\begin{example}\label{#1}}
\newcommand{\axe}{\end{example}}
\newcommand{\alg}{\begin{algorithm}}
\newcommand{\gla}{\end{algorithm}}
\newcommand{\nte}{\begin{note}}
\newcommand{\ntel}[1]{\begin{note}\label{#1}}
\newcommand{\etn}{\end{note}}
\newcommand{\app}{\begin{application}}
\newcommand{\appl}[1]{\begin{application}\label{#1}}
\newcommand{\ppa}{\end{application}}
\newcommand{\mat}{\begin{matrix}}
\newcommand{\tam}{\end{matrix}}
\newcommand{\smm}{\begin{summary}}
\newcommand{\mms}{\end{summary}}

\newcommand{\teo}{\begin{teorema}}
\newcommand{\teol}[1]{\begin{teorema}\label{#1}}
\newcommand{\oet}{\end{teorema}}
\newcommand{\df}{\begin{definicion}}
\newcommand{\dfl}[1]{\begin{definicion}\label{#1}}
\newcommand{\fd}{\end{definicion}}
\newcommand{\por}{\begin{proposicion}}
\newcommand{\porl}[1]{\begin{proposicion}\label{#1}}
\newcommand{\rop}{\end{proposicion}}
\newcommand{\cor}{\begin{corolario}}
\newcommand{\corl}[1]{\begin{corolario}\label{#1}}
\newcommand{\roc}{\end{corolario}}
\newcommand{\lem}{\begin{lema}}
\newcommand{\leml}[1]{\begin{lema}\label{#1}}
\newcommand{\mel}{\end{lema}}
\newcommand{\pru}{\begin{prueba}}
\newcommand{\prul}[1]{\begin{prueba}\label{#1}}
\newcommand{\urp}{\end{prueba}}
\newcommand{\axa}{\begin{axioma}}
\newcommand{\axal}[1]{\begin{axioma}\label{#1}}
\newcommand{\xa}{\end{axioma}}
\newcommand{\nta}{\begin{nota}}
\newcommand{\ntal}[1]{\begin{nota}\label{#1}}
\newcommand{\atn}{\end{nota}}
\newcommand{\eje}{\begin{ejemplo}}
\newcommand{\ejel}[1]{\begin{ejemplo}\label{#1}}
\newcommand{\je}{\end{ejemplo}}
\newcommand{\api}{\begin{aplicacion}}
\newcommand{\apil}[1]{\begin{aplicacion}\label{#1}}
\newcommand{\ipa}{\end{aplicacion}}

\newcommand{\ali}{\begin{align}}
\newcommand{\ila}{\end{align}}
\newcommand{\enu}{\begin{enumerate}}
\newcommand{\une}{\end{enumerate}}
\newcommand{\arr}{\begin{array}}
\newcommand{\rra}{\end{array}}
\newcommand{\eqa}{\begin{eqnarray}}
\newcommand{\aqe}{\end{eqnarray}}
\newcommand{\equ}{\begin{equation}}
\newcommand{\uqe}{\end{equation}}
\newcommand{\subq}{\begin{subequations}}
\newcommand{\qbus}{\end{subequations}}

\newtheorem{theorem}{Theorem}[section]
\newtheorem{lemma}[theorem]{Lemma}
\newtheorem{proposition}[theorem]{Proposition}
\newtheorem{corollary}[theorem]{Corollary}
\newtheorem{conjecture}[theorem]{Conjecture}

\theoremstyle{definition}
\newtheorem{definition}{Definition}[section]
\newtheorem{axiom}[definition]{Axiom}

\theoremstyle{remark}
\newtheorem{remark}{Remark}[section]
\newtheorem{example}[remark]{Example}
\newtheorem{application}[remark]{Application}
\newtheorem{algorithm}[remark]{Algorithm}
\newtheorem{note}[remark]{Note}
\newtheorem{summary}[remark]{Summary}

\newtheorem{teorema}[theorem]{Teorema}
\newtheorem{lema}[theorem]{Lema}
\newtheorem{proposicion}[theorem]{Proposici\'on}
\newtheorem{corolario}[theorem]{Corolario}

\newtheorem{definicion}[definition]{Definici\'on}
\newtheorem{axioma}[definition]{Axioma}

\newtheorem{nota}[remark]{Nota}
\newtheorem{ejemplo}[remark]{Ejemplo}
\newtheorem{aplicacion}[remark]{Aplicaci\'on}

\title{Experimental Evidence on a Refined Conjecture of the BSD type} 
\author{Francisco X. Portillo-Bobadilla}
\affil{Universidad Aut\'onoma de la Ciudad de M\'exico, Calzada Ermita Iztapalapa  No. 4163, Col. Lomas de Zaragoza, Iztapalapa, M\'exico  D.F. , C.P. 09620. \email{francisco.portillo@uacm.edu.mx}}

\begin{document}

\maketitle

\begin{abstract}
Let $E/\Q$ be an elliptic curve of level $N$ and rank equal to $1$. Let
$p$ be a prime of ordinary reduction. We experimentally study conjecture $4$ 
of B. Mazur and J. Tate in his article \emph{Refined Conjectures of the Birch and
Swinnerton-Dyer Type} \cite{bmjt}. 
We report the computational evidence.

\end{abstract}



\section{Introduction} \label{sect1}
B. Mazur and J. Tate in \emph{Refined Conjectures of the Birch and
Swinnerton-Dyer Type} postulated a series of conjectures of the
BSD-type in terms of finite layers. The goal was to find ``functions
with adelic type domains of definition and ranges of values'' for
which the $p$-adic $L$ functions were only a component, as expressed by
Yuri Manin \cite{man2}.
The Mazur and Tate conjecture (MT conjecture) is similar in spirit to the Birch and Swinnerton-Dyer conjecture (BSD conjecture). The conjecture has two assertion:
\enu
\item[1.] One that relates the rank of the elliptic curve with the order of vanishing of modular elements.  
\item[2.] The other that gives an explicit formula that relates arithmetic invariants of the curve with the modular element modulo the $r$-power of an augmentation ideal. In this formula, we have:
\enu
\item On the Arithmetic side: invariants like the Tamawaga constant, the order of the torsion group, the order of the Tate-Shafarevich group as exponents of a bi-multiplicative function, called the {\it corrected regulator}.

\item On the Analytic side: the modular element, defined in terms of modular symbols, and which is an analogue of a Stickelberger element.
\une

\une

In the present work, we show computational evidence only related to
the second assertion of the conjecture.
 
Our goal was to expand the evidence in favor of the conjecture (4) given by B. Mazur and J. Tate in \cite{bmjt}. In particular, they tested the conjecture for the elliptic curves $37$A and $43$A of rank 1 over sets $S=\{q\}$, where $q$ is a single prime of non split multiplicative reduction. They gave a very specific formula on those examples with prime conductor and group of Tate-Shafarevich trivial. We modify their equation so that any elliptic curve of rank 1 can be tested with no restrictions. 

The change consist on introducing adequate exponents on each side of the equation, the exponents depend on invariants of BSD type as the mentioned above, and we also introduce a value $\mu$ which is explained below. Hence, our contribution is to present a very concrete and easy to test conjecture and some computational evidence for it.

\section{Mazur-Tate Conjecture (General Setting)}\label{section-mt-conjecture}

Assume $E$ is an elliptic curve over $\Q$ with conductor $N$. Consider a N\'eron differential $\omega$ for $E$. Such $\omega$ is unique up to sign. 
Let $\Lambda_E$ be the N\'eron lattice (i.e.  
the lattice generated by the ``periods''
$\int_\gamma\omega\in\C$, where $\gamma$ runs through loops in
$E(\C)$) .

There is a unique pair of positive real numbers $\Omega_E^+$ and
$\Omega_E^-$ such that one of the two conditions holds:
\enu
\item $\Lambda_E=\Omega_E^+\Z+\Omega_E^-i\Z$
\item $\Lambda_E\subset\Omega_E^+\Z+\Omega_E^-i\Z$ is the sub-lattice
  generated by the complex numbers $a\Omega_E^+ +b\Omega_E^-i$ such
  that $a-b\equiv 0$ (mod $2$).
\une
In the first case, we say that $\Lambda_E$ is
  rectangular, otherwise $\Lambda_E$ is non-rectangular.
 
Let $f$ be the modular form associated to $E$, and let $a/b$ be a
rational number. We define the modular elements $[a/b]^+_E$ and $[a/b]^-_E$ by:
\equ \label{eq1}
2\pi\int_{0}^{\infty}f(a/b+it)dt=\Omega_E^+[a/b]^+_E+\Omega_E^-[a/b]^-_E
i.
\uqe

We will write $[a/b]$ instead of $[a/b]^+_E$, since we will be concerned only with the plus symbols on $E$. The number $[a/b]$ is rational, and if $b$ is prime to the conductor of the curve, the value is an integer \cite{man1}.

Let $S$ be a finite set of primes, let $S'$ be the subset of $S$ of primes with multiplicative reduction at a fixed elliptic curve $E$. Set
\equ
M=\prod_{p\in S-S'}p\prod_{p\in S'}p^{e_p}
\uqe
with integers $e_p\geq 0$, and set
\equ
G_M=(\Z/M\Z)^*/(\pm 1)
\uqe
If $a$ is an integer coprime to $M$, let $\sigma_a$ denote its
associated element in $G_M$.

Let $R$ be a subring of $\Q$ containing $1/2$ and $1$ over the order of the torsion of $E(\Q)$. Define the modular element as
\equ
\Theta_{E,M}:=\frac{1}{2}\sum_{a \mathrm{mod} M}\left[\frac{a}{M}\right]\cdot\sigma_a\in R[G_M]
\uqe

Let $\epsilon:R[G_M]\rightarrow R$ be the augmentation map,
defined by 
\equ
\sum r_i \sigma_a \rightarrow \sum r_i
\uqe
and let $I=\ker(\epsilon)$ its augmentation ideal.

Let $X$ be the N\'eron model of $E$, let $X(\F_p)$ be fiber of the N\'eron model of $E$ at $p$,
let $X^0(\F_p)=E_{ns}(\F_p)$ be the non-singular points of $E$ modulo $p$ and
let $N_p=X(\F_p)/X^0(\F_p)$ be
the group of connected components in the fiber.

Define $\phi_{S}$ as the order of the cokernel of the natural projection:
\equ
\pi_{S} : E\rightarrow\prod_{p\notin S'}N_p
\uqe
as $q$ ranges through the set of all primes.

Conjecture $4$ in \cite{bmjt} is the following:

\cnj \label{mtconj}(``Birch-Swinnerton-Dyer type'' conjecture.)

Let $r=\mathrm{rank}(E(\Q))+\#(S')$. The modular element $\Theta_{E,M}\in R[G_M]$ lies in
the $r$-th power of the augmentation, $I^r\subset R[G_M]$, and if $\tilde{\Theta}_{E,M}$ denotes its image in $I^r/I^{r+1}$:
\equ
\tilde{\Theta}_{E,M}=\#(\Sh)\cdot\phi_{S}\cdot\nu_r(\mathrm{Disc}_S(E))\in I^r/I^{r+1} 
\uqe
\jnc 

In the following pages, we explain the term $\nu_r(\mathrm{Disc}_S(E))$.

\subsection{Definition of $\mathrm{Disc}_S(E)$.}

\subsubsection{Local construction of the regulator}\label{local_regulator}





Using the theory of biextensions and splittings, Mazur and Tate introduce local canonical heights and corrected discriminants. We give a brief summary of their work to introduce regulators. For more details, see \cite{bmjt0} and \cite{bmjt}.


\begin{definition}
If $A$, $B$ and $C$ are abelian groups. A biextension of $(A,B)$ by $C$ is an object $\Ec$ such that for each triple
$(a,b,c)\in A\times B\times C$, we can assign a unique element $[a,b,c]\in \Ec$ such that $a\Ec:=[a,B,C]\subseteq \Ec$ has a group structure isomorphic to $B\times C$; and analogously, $b\Ec:=[A,b,C]$ has a group structure isomorphic to $A\times C$. Also, $C$ acts freely on $\Ec$. 
\end{definition}

Now, let $\tilde{A}$, $\tilde{B}$ and $\tilde{C}$ be other abelian groups. If $\alpha:\tilde{A}\rightarrow A$, $\beta:\tilde{B}\rightarrow B$ are
injective homorphisms, and  $\rho: C\rightarrow \tilde{C}$ is a surjective homomorphism, we can obtain a biextension $\tilde{\Ec}$ given by the pullback of $\Ec$ by $\alpha$ and $\beta$, and the pushout of $\Ec$ by $\rho$. 

\begin{definition}
Let $\Ec$ be a biextension of $(A,B)$ by $C$, and $\rho:C\rightarrow \tilde{C}$ a group homomorphism. A $\rho$-splitting of $\Ec$ is a map $$\psi:\Ec\rightarrow \Ct$$ such that
\enu
\item $\psi(\omega\cdot x)=\rho(\omega)\cdot \psi(x)$ for $x\Ec$ and $w\in C$.
\item $\psi|_{a\Ec}$ and $\psi|_{b\Ec}$ are group homomorphisms.
\une
\end{definition}

If $A$ and $B$ are dual varieties over a field $K$, we know that there exists a bi-extension $\Ec$ of $(A,B)$ by $K^*$
that expresses the duality \cite{gr7}. Denote this biextension by $\Ec(K)$.

\begin{definition}
A modification $(\tilde{\Ec},\alpha,\beta,\rho)$ of $\Ec(K)$ is a biextension $\tilde{\Ec}$ obtained by injective homomorphisms $\alpha:\tilde{A}\rightarrow A(K)$, $\beta:\tilde{B}\rightarrow B(K)$; and a surjective homomorphism  $\rho: K^*\rightarrow \tilde{C}$.
\end{definition}

\begin{definition}
A trivialization $(\alpha,\beta,\rho,\psi)$ of $\Ec(K)$ is a modification
$(\tilde{\Ec},\alpha,\beta,\rho)$ of $\Ec(K)$ and a $\rho$-splitting $\psi$ of $\tilde{\Ec}$.
\end{definition}

Notice that if $\alpha:\tilde{A}\rightarrow A(K)$, $\beta:\tilde{B}\rightarrow B(K)$ and $\rho:K^*\rightarrow C$ are group homomorphism as above, and $(\tilde{\Ec},\alpha,\beta,\rho)$ is the associated modification, we have a bi-multiplicative function
$$\langle\mbox{ },\mbox{ }\rangle_{\tilde{\Ec}}:\tilde{A}\times\tilde{B}\rightarrow \tilde{C}$$ defined by
$$\langle\tilde{a},\tilde{b}\rangle_{\tilde{\Ec}}:=\rho\left(\langle\alpha(\tilde{a}),\beta(\tilde{b})\rangle_\Ec\right) \mbox{ ($\forall\tilde{a}\in\tilde{A}$ and $\forall\tilde{b}\in\tilde{B}$) }.$$
Here, 
$\langle\mbox{ },\mbox{ }\rangle_\Ec:A\times B\rightarrow K^*$ is the bilinear pairing that express the duality.



If we define $\tilde{\psi}: \tilde{\Ec}\rightarrow \tilde{C}$ as $$\tilde{\psi}\left(\left[\tilde{a},\tilde{b},\tilde{c}\right]\right):=\tilde{c}\cdot\langle\tilde{a},\tilde{b}\rangle_{\tilde{\Ec}},$$
thus $\tilde{\psi}$ is a $\rho$-splitting of $\tilde{\Ec}$.
And therefore, $(\tilde{\Ec},\alpha,\beta,\rho)$  is a trivilization of $\Ec(K)$.

Working over local fields, Mazur and Tate \cite{bmjt} described what they called
{\it ``the canonical trivilizations"}. From now on, we will assume that our local fields are the fields $\Q_p$ for $p$ a prime number, that our global field is $K=\Q$, that $A=E$ is an elliptic curve and $B=E^\vee$ is its dual variety.
Also, for each prime $p$, we will consider a system of group homomorphisms:
$$\alpha_p:A_p\rightarrow E(\Q_p),$$ $$\beta_p:B_p\rightarrow E(\Q_p),$$ 
$$\rho_p:\Q_p^*\rightarrow C_p,$$
where $\alpha_p$ and $\beta_p$ are injective and $\rho_p$ is surjective.

Hence, we will have modifications $(\alpha_p,\beta_p,\rho_p)$ with their corresponding $\rho_p$-splittings. For the purpose of this article, we are interested in the following three
trivializations:
\enu
\item[a)] {\it N\'eron unramified trivialization.} 
Let $A_p=E(\Q_p)$, $B_p=E^0(\Q_p)$ and $C_p=\Q_p^*/\Z_p^*\simeq \Z$.
Here, $E^0(\Q_p)$ denotes the group of points in $E(\Z_p)$ whose reduction modulo
$p$ is in the componente of zero in the fiber $E(\F_p)$.
The homomorphisms $\alpha_p$ and $\beta_p$ are the natural inclusions;
$\rho_p:\Q_p^*\rightarrow \Q_p^*/\Z_p^*$ is the natural projection.
Now, $\psi_p:\Ec_p\rightarrow \Q_p^*/\Z_p^*$ is the only canonical splitting
such that $\psi(\Ec_p(\Z_p))=0$. 

\item[b)] {\it Tamely ramified trivialization.}
Let $A_p=E(\Q_p)$, $B_p=E^1(\Q_p)$, $C_p=\Q_p^*/p\Z_p^*\simeq \F_p^*$.
The maps $\alpha_p$ and $\beta_p$ are the inclusions again and
$\rho_p$ is the projection.
Now, $E^1(\Q_p)$  are the points in $E(\Q_p)$ whose reduction modulo $p$ is zero in the 
conected component of zero in the fiber $E(\F_p)$.

\item[c)] {\it Split Multiplicative trivialization.}
If $p$ is a prime of split multiplicative reduction, then $E(\Q_p)$ is isomorphic to the Tate curve $E_{q_p}=\Q_p^*/q_p^\Z$, where $q_p$ is the multiplicative local period. Hence, in this trivialization, 
we take $A_p=B_p=\Q_p^*$ and $C_p=\Q_p^*$. And, $\beta_p=\alpha_p:\Q_p^*\rightarrow\Q_p^*/q_p^\Z$ is the natural parametrization of $E_{q_p}$, and 
$\rho_p:\Q_p^*\rightarrow C_p=\Q_p^*$ is the identity.
\une


\subsubsection{Global Construction of Regulator}

For the finite set $S$ (See section \ref{section-mt-conjecture}.), we will construct extended Mordel groups
$A_S$, $B_S$ and $C_S$ as follows:

According to subsection \ref{local_regulator}, for each subset of primes $S\subseteq\wp$, there is a system of homomorphisms
$$\alpha_p:A_p\rightarrow E(\Q_p),$$ $$\beta_p:B_p\rightarrow E(\Q_p),$$ 
$$\rho_p:\Q_p^*\rightarrow C_p,$$ with their corresponding trivializations
$$\psi_p:\tilde{\Ec}_p\rightarrow C_p.$$ 

The trivialization $\psi_p$ is determined by the rule:

\enu
\item[a)] $\psi_p$ is the {\it N\'eron unramified trivialization}, if $s\notin S$.

\item[b)] $\psi_p$ is the {\it Tamely ramified trivialization}, if $s\in S-S'$.

\item[c)] $\psi_p$ is the {\it Split Multiplicative trivialization}, if $s\in S'$.

\une

We define $A_S$ to be the set of pairs $(P,(a_p))$ such $P\in E(\Q)$, $(a_p)\in\prod_{p\in\wp}A_p$ and $\alpha_p(a_p)=i_p(P)$ for all prime $p$, where 
$i_p:A(\Q)\rightarrow E(\Q_p)$ is the canonical inclusion. We define $B_S$, similarly.
 

Now, from the $3$ possibilities of local trivializations, we can write
$C_p=\Q_p^*/U_p$, where $U_p$ could be either $\Z_p^*$, $p\Z_p^*$ or $\{1\}$.

Hence, we have a morphism $$\rho:=(\rho_p):\coprod_{p\in\wp}\Q_p^*\rightarrow\bigoplus_{p\in\wp}(\Q^*_p/U_p).$$

Now, if we mod out by $\Q^*$ using the natural inclusions $\Q^*\hookrightarrow\Q_p^*$, define: 
\equ
C_S:=\coprod_{p\in\wp}\Q_p^*/\Q^*(\prod_{p\in\wp}U_p).
\uqe 
Set
$$\phi:\bigoplus_{p\in\wp} C_p\rightarrow C_{S},$$
the natural map given by coordinates.

For $a=(P,(a_p))\in A_S$ and $b=(Q,(b_p))\in B_S$, 
define the bimultiplicative
pairing by
\equ
\langle a,b\rangle_S:=\phi(\prod_{p}\psi_p(x_p))=\prod_{p}\phi\circ\psi_p(x_p),
\uqe
where $x_p=[a_p,b_p,k]\in\tilde{\Ec}_p$

Notice 
$A_p= E(\Q_p)$
and $B_p= E^0(\Q_p)$ for almost all $p$, and since $P\in A(\Z_p)$ and $Q\in B(\Z_p)$ for almost all $p$, we have that $\psi_p(x_p)=1$ for almost all $p$.

Hence, the global bi-multiplicative function 
is computed as a finite product.

In fact, in our example, we have
\equ\label{C_S_1}
C_S:=\left(\prod_{p\in S-S'}\F_p^*\times \prod_{p\in S'}\Z_p^*\right)/(\pm 1) .
\uqe 

Now, $A_S$ and $B_S$ are finitely generated groups of the same rank:
$$r=\mathrm{rank}(A(K))+\#(S')\cdot \dim(A)\mbox{ .  (Reference: \cite{bmjt}.)}$$

Hence, if $\{P_1,P_2,\ldots,P_r\}$ generates the free part of $A_S$ and $\{Q_1,Q_2,\ldots,Q_r\}$ generates the free part of $B_S$, set
$$\mathrm{disc}_S=\det_{1\leq i, j\leq r} \langle P_i, Q_j \rangle$$

The value $\mathrm{disc}_S$ is well defined up to sign. But, we can choose an adecuate
orientation for our purposes.

Now, for our computations, it is useful to work on a subring $R\subset\Q$ containing the torsion of $A_S$ and $B_S$.  Hence, we will consider the element
$$d_S:=1\otimes \mathrm{disc}_S \in R\otimes \mathrm{Sym}_r(C_S).$$

This discriminant does not work well as the regulator, see the heuristic discussion about it in \cite{bmjt}.

Instead, the {\it corrected discriminant} is defined as a sum of discriminants
$d_T$ over subsets $T\subset S$ containing $S'$.

For any subset $T\subset S$, we have natural mappings:
$x_{S,T}:A_S\rightarrow A_T$, $y_{S,T}:B_S\rightarrow B_T$ and
$z_{T,S}:C_T\rightarrow C_S$. 

There is also a unique map $\mu_{S,T}:C_S\rightarrow C_T$, such that
$$\mu_{S,T}\circ z_{T,S}=\prod_{p\in S-T} (p-1)\cdot c\mbox{ }\mbox{ for all $c\in C_S$}$$

Thus, the {\it corrected discriminant} of $S$ is defined as:

\equ\label{corrected_disc}
\mathrm{Disc}_S(A)=\sum_{S'\subset T\subset S} (-1)^{\#(T-S')} \mu_{S,T}(j_T\cdot d_T)\in R\otimes \mathrm{Sym}_r(C_S),
\uqe
where 
$j_T=\left(\prod_{p\in S-S'} n_p \right)/ (B_{S'}:B_S)$, $n_p=\# B^0(\Q_p)$, and
$(B_{S'}:B_S)$ is the index of $B_S$ in $B_S'$.

Now, from equation (\ref{C_S_1}) there is a natural surjective homomorphism $C_S \twoheadrightarrow G_M$. And, also a natural identification of $G_M$ with
$I^2/I$ (as is described in next section).
Thus, we have a natural map $C_S\rightarrow I^2/I$, which induces a natural
homomorphism:
\equ
\nu_r: R\otimes \mathrm{Sym}_r(C_S)\rightarrow I^{r+1}/I^r
\uqe
 
Now, we should notice that the formula in Conjecture \ref{mtconj} ocurrs
in $I^{r+1}/I^r$, and thus, the analogous of the regulator is $\nu_r\left(\mathrm{Disc}_S(A)\right)$.

\section{MT Conjecture (Rank 1, Ordinary and Good Reduction Setting)}

\subsection{The Analytic Side}

In this section, we assume that $E$ has rank $1$ and that $S$ has only primes of ordinary reduction. In this context, Conjecture \ref{mtconj} in section \ref{section-mt-conjecture} states
that
\enu
\item[a)] $$\Theta_{E,M}\in I$$
\item[b)] $$\tilde{\Theta}_{E,M}=\#(\Sh)\cdot\phi_{S_m}\cdot\nu_r(\mathrm{Disc}_S(E))\in I/I^{2}$$
\une

Now, assertion a) is equivalent to have $\epsilon(\Theta_{E,M})=0$, or equivalently
\equ
\sum_{a \mathrm{mod} M}\left[\frac{a}{M}\right]=0.
\uqe

Hence, we have
\equ
\Theta_{E,M}=\frac{1}{2}\sum_{a \mathrm{mod} M}\left[\frac{a}{M}\right]\cdot(\sigma_a-e)\in R[G_M]
\uqe
where $e$ is the identity on $G_M$.

The Hurewicz Theorem for augmentation ideals gives an isomorphism of abelian groups
$G_M\simeq I/I^2$ given by
the map $r(g-e)\mapsto g^r$ for $g\in G$ and $r\in \Z$. Hence, we will test assertion b) of the Conjecture directly on the group:
$$G_M=\left(\prod_{p\in S}\F_p^* /\pm 1\right).$$

Since we cannot compute always square roots in $\F_p^*$, we will test the conjecture for the square of $\tilde{\Theta}_{E,M}$, which is equivalent to eliminate the $\frac{1}{2}$ on $\Theta_{E,M}$.
Conjecture \ref{mtconj} is additive, but our testing will be multiplicative.

\begin{definition}\label{def_mod_symb}
For $S$ having only primes of good reduction and an elliptic curve $E$ with $\mathrm{rank}(E)\geq 1$, we define the following multiplicative modular
element:
\equ \label{eq4}
l(S)=\prod_{a\in (Z/MZ)^*} a^{[a/M]} \mbox{(mod $M$)}
\uqe
\nfd
with $M=\prod_{p\in S} p$ .

The values $[a/M]$ are integers if $gcd(M, N)=1$ by 5.4 in  \cite{man1}, so the
multiplicative modular element is well defined.



\subsection{The Arithmetic side}

In this section, we also assume that $E$ is an elliptic curve with positive rank.
First, assume $p$ is a prime of good reduction and $S=\{p\}$. In this case, we will describe how to compute $\mathrm{Disc}_S(E)$.

An element $x\in\Ec_p$, can also be described by a triplet
$x=[\mathfrak{a},D,c]$, where $\mathfrak{a}=\sum_i n_i(P_i)$ is a zero
cycle with $P_i\in E(\Q_p)$, $D=\sum_j m_j(Q_j)$ is a divisor in $E^0(\Q_p)$ algebraically  equivalent to zero whose support is disjoint to $\mathfrak{a}$, and $c\in\Q_p^*$ \cite{bmjt} and \cite{neron}.

Now, 
this symbol satisfies the properties:
\enu
\item[a)] $[\mathfrak{a},div(f),1]=[\mathfrak{a},0,f(\mathfrak{a})]$ for a rational function
$f$ defined on $E(\Q_p)$ with $f(\mathfrak{a})=\prod_jf(Q_j)^{m_j}$.
\item[b)] $[\mathfrak{a}_R,D_R,c]=[\mathfrak{a},D,c]$, where $\mathfrak{a}_R$ (resp. $D_R$) is obtained from $\mathfrak{a}$ (resp. $D$) by translating each point by $R$.
\une

Now, since $E$ is an elliptic curve, we identify a point $P\in E$, with the zero cycle
$(P)-(O)$. Hence, the discriminant is
\equ
\mathrm{Disc}_{\{p\}}(E)=\psi_p([(P)-(O),(O)-(Q_p),1])
\uqe
where $P$ is a generator of $E(\Q)$, $Q$ is a generator of $E^0(\Q)$ and 
$Q_p=n_p Q$.

Notice that the element $[(P)-(O),(O)-(Q_p),1]$ is above $E(\Q_p)\times E^1(\Q_p)$
on the biextension $\Ec_p$.

Now, the value $\psi_p([\mathfrak{a},D,1])$ coincides with the N\'eron's symbol
$(D,\mathfrak{a})_{v_p} $ for the $p$-adic valuation in $\Q_p$. In particular, Theorem $3$ in \cite{neron} says
how to compute  $(D,\mathfrak{a})_{v_p} $ if $D$ is equivalent to $O$.

To compute $\mathrm{Disc}_{\{p\}}(E)$ is helpful to use property 2) above, translating by a point $P'$. Hence,
\equ
\mathrm{Disc}_{\{p\}}(E)=\psi_p([(P+P')-(P'),(P')-(Q_p+P'),1])
\uqe
 
This value is the $g$ function defined by Mazur and Tate in page $747$ of \cite{bmjt}:

Let $P$, $Q$ and $P'$ be as above. For  $p\nmid N$ prime, consider the
quantity:

\equ \label{eq6}
g(P,Q,P',p)=\frac{d(P'+P)d(P'+Q_p)}{d(P')d(P'+P+Q_p)}\mbox{ (mod $p$)}
\uqe
where $d(T)$ is the square root of the denominator of the $x$-coordinate of a point $T$.

We will consider the square of this $g$ function, just assuming that $d(T)$ is the
$x$-coordinate of $T$. This will balance the cancellation of the $\frac{1}{2}$
in $\Theta_{E,M}$, and it is in concordance with definition \ref{def_mod_symb}.


We sumarize the properties of the $g$ function in the following proposition:


\pro \label{mrlema}\mbox{ }

\enu
\item If $P\in E(\Q)$, $Q\in E^0(\Q)$, 
then $g(P,Q,P',p)$ does not depend on $P'$. Moreover, if $P$ is a generator of the free part of $E(\Q)$ and $Q$ is a generator of the free part of $E^0(\Q)$, then this value depends only on $E$ and $p$.
\item The function
$$
\hat{g}:E\times E_0\rightarrow \prod_{p\nmid N} \F_p^*,
$$
defined by $$\hat{g}(P,Q)_p:=g(P,Q,P',p)\mbox{ at the $p$ coordinate}$$ is bi-multiplicative.
\une
\orp


Now, let $S$ be a finite set of primes having only good reduction at $E$.
Set $M=\prod_{p\in S}p$, $n_S=\prod_{p\in S}n_p$ and $Q_S=n_S Q$.
If $P$ and $Q$ are generators of the free part of $E$, define 
\equ \label{gfunction}
g(S)=\frac{d(P'+P)d(P'+Q_S)}{d(P')d(P'+P+Q_S)}\mbox{ (mod $M$)}
\uqe
where $P'$ is a point on $E$ such than non of the $d$'s is zero.

Now, if $M'\mid M$, let 
$$Y_{M',M}:(\Z/M'\Z)^*\rightarrow(\Z/M\Z)^*$$
be the map defined by  $a\mapsto b^{\phi(M/M')}$,
where
$a\in(\Z/M'\Z)^*$ and $b\in(\Z/M\Z)^*$ such that $a\equiv b$ (mod $M'$), and
$\phi$ is {\it the Euler phi}.

\begin{definition}
The $G$ function on $S$ is
\equ
G(S):=\prod_{T\subseteq S} Y_{M_T,M}\big(g(P,Q_T,M_T)\big)^{(-1)^{(1+\#(T))}}
\uqe
where $M_T=\prod_{q\in T}q$, $n_T=\prod_{p\in T}n_p$ and $Q_T= n_T Q$.
\end{definition}


\subsection{Multiplicative Equations of Mazur-Tate Conjecture}

Assume $E$ is an elliptic curve of rank $1$. Let $E_0$ be the group of everywhere good reduction
points of $E$. 

First, assume $S$ has only points of ordinary reduction (i.e. $S'=\{\}$). Therefore, $\phi_S$ 
is the cokernel of the natural projection:
\equ
\pi_{S}: E\rightarrow\prod_{p\in\wp}N_p
\uqe
where $p$ ranges through the set of all primes $\wp$.

The kernel of $\pi_S$ is $E_0$.
Hence, the induced map $$E/E_0\hookrightarrow \prod_{p\in\wp}N_p$$ is an injection of finite groups and its cokernel is the cokernel of $\pi_S$. Hence,




\equ
\phi_S=\frac{C}{\#(E/E_0)},
\uqe
where $C=\#\left(\prod_{p\in\wp}N_p\right)=\prod_{p\in\wp}c_p$ and
$c_p=|N_p|$ are the Tamagawa numbers.

If $S'\neq \emptyset $, then we divide by $C'=\prod_{p\in S'} c_p$, to obtain  
\equ\phi_S=\frac{C}{C' \#(E/E_0)}.\uqe

Let $E_{tors}$ be the group of torsion points of $E$. 
If $u$ is the order of torsion in $E$ and $v$ is the order of the torsion in $E_0$,
then we can explicitly compute the order $\#(E/E_0)$ as
$\frac{\mu u}{v}$, where
\equ
\mu=\mbox{min\{$j>0$ : $jP+R\in E_0$ and $R\in E_{tors}$\}}
\uqe
and $P$ is any generator of the free part of $E$.

Thus, Conjeture \ref{mtconj} on its multiplicative form and running over all good reduction points gives:

\cnj \label{mazur2} (Rank $1$ at all Good Reduction Primes.)

Let $E$ be a curve of rank $1$, let $P$ be a generator of $E$ (modulo
torsion), and let $Q$ be a generator of $E_0$ (modulo torsion), then:
\begin{equation}
\hat{l}^{uv}=\hat{g}(P,Q)^{|\Sh| |coker(\phi)|}\in\prod_{p\nmid N} \F_p
\end{equation}
where $|\Sh|$ is the order of the Tate-Shafarevich group and
$\hat{l}=\prod_{p\nmid N} l(\{p\})$.

\jnc


Notice that if we exponentiate the above equation by $u/v$, we obtain the
equation:
\equ
\hat{l}^{u^2}=\hat{g}(P,Q)^{\frac{C |\Sh |}{\mu}}
\uqe
which looks more like the classical BSD.

For a more general $S$, having only good reduction points, the conjecture \ref{mtconj} in its multiplicative form becomes:

\cnj \label{mazur3} (Rank $1$ for S having only Good Reduction Primes.)

Let $E$ be a curve of rank $1$ and $S$ having only good reduction primes, then:
\begin{equation}
l(S)^{uv}=G(S)^{|\Sh| |coker(\phi)|}\in G_M
\end{equation}
where $M=\prod_{p\in S}p$ and $|\Sh|$ is the order of the Tate-Shafarevich group.

\jnc

In Chapter 4 of \cite{port}, we explained how to test Conjecture \ref{mazur3} using the individual computations
on each prime $p\in S$.


\section{Testing conjectures \ref{mazur2} and \ref{mazur3}.}


On \cite{port}, we tested the above conjecture for the first $300$ elliptic curves in
the Cremona database \cite{cre1}. All these cases have trivial
Tate-Shafarevich group. But, we also tested in \cite{port} for an elliptic curve 
having a non-trivial
Tate-Shafarevich group. The curve was
\equ
y^2+xy+y=x^3-x^2-8587x-304111
\uqe
with conductor $N=1610$ and $|\Sh|=4$. 

Those computations were done using the Pari calculator \cite{pari} with the help of
the script \cite{msym}, we tested each curve for $p<300$ and $p\nmid N$.

Now, we enlarge our experimental evidence using
SAGE \cite{sage}. We test the Conjecture \ref{mazur2} on the
first $3000$ curves elliptic on the Cremona database (already included in SAGE). 

We also check the Conjecture \ref{mazur2} for more elliptic curves with non-trivial Tate-Shafarevich group. We check on the first $20$ elliptic curves with $|\Sh|=4$ and on the first $7$ elliptic curves with $|\Sh|=9$. We use
{\it The L-functions and Modular Forms Database} \cite{mfdb} to search for the required elliptic curves to test.  

The files with the computing evidence and the scripts are available on

\mbox{ }

\framebox[3.3in]{ https://github.com/portillofco/MazurTateProject }

\mbox{ }

\nte {\bf Last comment regarding normalization of modular symbols.}

We use the usual methods for computing modular symbols and take advantage of the computing power of Pari-gp and Sage. There have been continous advancement on the methods for computing modular symbols and also in the computing power used on computations, but correct normalization is still a practical issue to be considered during the testing of the conjecture. 

The computation of the modular symbols $[a/b]^+$ using only Linear Algebra is alright up to multiplication by a constant. On our first computations \cite{port} using Pari, we determined the constant by a series aproximation of the value $[a/b]^+$. 
Now, Sage computes $[a/b]^+$ correctly in most of the cases, but there are still
a few curves when Sage prompts a WARNING MESSAGE. 

For example, for the curve 158 in the Cremona Data Base, we received the following WARNING MESSAGE:

{\bf Warning : Could not normalize the modular symbols, maybe all further
results will be multiplied by -1, 2 or -2.}

In such cases, we just verified which of the proposed values works for the conjecture. We must point out that in all the curves tested, one of the suggested values works. We believe that some numerical modular symbols can be used to compute the constant in a direct way \cite{wuth}.

\etn

Finally, we mention that we made the computations using a HP Workstation with a Procesor Intel Xeon E5-2640v2 with 8 nodes and 48GB of RAM memory.










\subsubsection*{Acknowledgements}
I would like to thank Felipe Voloch for his guidance and advice during the development of this research. I am also very grateful
to John Tate for his valuable help and his many explanations regarding
the conjecture $4$ in \cite{bmjt}. The updating of this article and the new experimental
evidence was supported by the project {\it PI2013-38} of the agreement 
{\it UACM/SECITI/060/2013}. I thank also the support of my collegues Isa\'{\i}as L\'opez and Felipe Alfaro during the development of the aforementioned project.

\end{document}